\documentclass[makeidx]{article}
\usepackage{makeidx}
\usepackage{amssymb , amsmath, amsthm}
\newtheorem{defi}{Definition} 
\newtheorem{thm}[defi]{Theorem}

\newtheorem{rem}[defi]{Remark}
\newtheorem{prop}[defi]{Proposition}
\newtheorem{lemme}[defi]{Lemma}
\newtheorem{cor}[defi]{Corollary}
 
\newcommand{\twosystem}[2]{\left\{\begin{aligned} &#1\\ &#2\end{aligned}\right.}

\newcommand{\nero}{\smallskip$\bullet\quad$\rm}

\newcommand{\parte}[1]{\smallskip\noindent {\rm#1)}\,\,}
\newcommand{\acapo}{\smallskip\noindent}

\newcommand{\scal}[2]{\langle{#1},{#2}\rangle}
\newcommand{\starred}[1]{{#1}^{\star}}

\newcommand{\Cal }[1]{{\mathcal {#1}}}
\newcommand{\bd}{\partial}

\newcommand{\norm}[1]{\lVert{#1}\rVert}

\newcommand{\reals}{{\bf R}}
\newcommand{\real}[1]{{\bf R}^{#1}}
\newcommand{\sphere}[1]{{\bf S}^{#1}}

\newcommand{\Vol}{{\rm Vol}}

\begin{document}

\title{A Reilly formula and eigenvalue estimates for differential forms\footnote{Classification AMS $2000$: 58J50, 58J32, 53C24 \newline 
Keywords: Manifolds with boundary, Hodge Laplacian, Spectrum, Rigidity}} 
\author{ S. Raulot and A. Savo}
\date{\today}
\maketitle
 
\begin{abstract}
We derive a Reilly-type  formula for differential $p$-forms on a compact manifold with boundary and apply it  to give a sharp lower bound of the spectrum of the Hodge Laplacian acting on differential forms of an embedded hypersurface of a Riemannian manifold. The equality case of our inequality gives rise to a number of rigidity results, when the geometry of the boundary has special properties and the domain is non-negatively curved. Finally we also obtain, as a by-product of our calculations, an upper bound of the first eigenvalue of the Hodge Laplacian when the ambient manifold supports non-trivial parallel forms. 
\end{abstract}

\large


\section{Introduction}


Let $\Omega$ be a compact Riemannian manifold of dimension $n+1$ with smooth boundary $\Sigma$, and $f$ a smooth function on $\Omega$. In \cite{reilly},  Reilly integrates  the Bochner formula for the $1$-form $df$ and  re-writes the boundary terms in a clever way to obtain what is known  in the literature as the \it Reilly formula \rm for the function $f$.
Reilly used it to prove the Alexandrov theorem, but the Reilly formula turned out to be extremely useful, and has been successfully applied in many other contexts. For example, it was applied by Choi and Wang in \cite{choi-wang}, to obtain a lower bound for the first eigenvalue $\lambda_1(\Sigma)$ of the Laplacian on functions on embedded minimal hypersurfaces of the sphere, and by Xia in \cite{xia}, who proved the following lower bound: if the domain $\Omega$ has non-negative Ricci curvature and its boundary $\Sigma$ is convex, with principal curvatures bounded below by $c>0$, then:
\begin{equation}\label{xia}
\lambda_1(\Sigma)\geq nc^2;
\end{equation}
moreover the equality holds if and only if $\Omega$ is a ball of radius $\frac 1c$ in Euclidean space.   

In this paper we integrate the Bochner-Weitzenb\"ock formula to obtain a Reilly-type  formula for differential $p$-forms (see Theorem \ref{Reilly}) and we apply it to give a sharp lower bound of the spectrum of the Hodge Laplacian acting on differential forms of an embedded hypersurface $\Sigma$ bounding a compact manifold $\Omega$. This lower bound can be seen as a generalization of Xia's estimate to differential forms. 
The equality case gives a number of rigidity results.

The Reilly formula has been generalized in other contexts, for example in spin geometry (see \cite{hijazi-montiel-zhang}); we hope that the present paper fills a gap in the literature and that other applications could be found. 

Let us give a brief overview of the results of the paper.


\subsection{The main estimate}


Let us state the main lower bound in precise terms. Fix a point $x\in\Sigma$ and let $\eta_1(x),\dots,\eta_n(x)$ be the principal curvatures of $\Sigma$ with respect to the inner unit normal. The \it $p$-curvatures \rm are by definition all possible sums $\eta_{j_1}(x)+\dots+\eta_{j_p}(x)$ for increasing indices $j_1,\dots,j_p\in\{1,\dots,n\}$. Arrange the sequence so that it is non-decreasing:
$
\eta_1(x)\leq\dots\leq\eta_n(x);
$
then the \it lowest $p$-curvature \rm is $\sigma_p(x)=\eta_1(x)+\dots+\eta_p(x)$, and we set:
$$
\sigma_p(\Sigma)=\inf_{x\in\Sigma}\sigma_p(x).
$$
We say that $\Sigma$ is $p$-convex if $\sigma_p(\Sigma)\geq 0$, that is, if all $p$-curvatures are non-negative. Note that $1$-convex means convex (all principal curvatures are non-negative) and $n$-convex means that $\Sigma$ has non-negative mean curvature. It is easy to verify that, if $\sigma_p\geq 0$, then $\sigma_q\geq 0$ for all $q\geq p$ and moreover $\frac{\sigma_p}{p}\leq\frac{\sigma_q}{q}$. 

Here is the main estimate of this paper.  

\begin{thm}\label{Hodge}
Let $\Omega$ be a compact $(n+1)$-dimensional Riemannian manifold with smooth boundary $\Sigma$, and let $1\leq p\leq \frac{n+1}2$. Assume that $\Omega$ has non-negative curvature operator and the $p$-curvatures of $\Sigma$ are bounded below by $\sigma_p(\Sigma)>0$. Then:
$$
\lambda_{1,p}'(\Sigma)\geq
\sigma_p(\Sigma)\sigma_{n-p+1}(\Sigma),
$$
where $\lambda_{1,p}'(\Sigma)$ is the first eigenvalue of the Hodge Laplacian acting on exact $p$-forms of $\Sigma$. The equality holds if and only if $\Omega$ is a Euclidean ball. 
\end{thm} 

We will prove the lower bound under weaker curvature conditions, namely that the curvature term in the Bochner formula for $p$-forms on $\Omega$, denoted $W^{[p]}_{\Omega}$, is non-negative (we refer to section \ref{ReillyFormula} and Theorem \ref{Hodge2} for this more general result).

\smallskip

By the Hodge duality one has $\lambda_{1,p}'(\Sigma)= \lambda_{1,n-p+1}'(\Sigma)$, so the positivity of $\sigma_p(\Sigma)$ for some $1\leq p\leq \frac{n+1}2$ allows to bound from below the eigenvalues $\lambda_{1,q}'(\Sigma)$ for all $q\in [p,n-p+1]$. In particular, if $\Sigma$ is convex, with principal curvatures bounded below by $c>0$, then it is $p$-convex for all $p$ and $\sigma_p(\Sigma)\geq pc$. Therefore, for all $p=1,\dots,n$ one has:
$$
\lambda_{1,p}'(\Sigma)\geq p(n-p+1)c^2.
$$
For $p=1$, we get back Xia's estimate \eqref{xia} because $\lambda_{1,1}'(\Sigma)=\lambda_1(\Sigma)$, the first positive eigenvalue of the Laplacian acting on functions.  
\medskip

\begin{rem}\label{wu}\rm
The $p$-convexity has interesting topological consequences: Wu proved in \cite{wu} that if $\Omega$ has non-negative sectional 
curvature and $\Sigma$ is strictly $p$-convex (that is, $\sigma_p(\Sigma)>0$), then $\Omega$ has the homotopy type of a CW-complex with cells only in dimension $q\leq p-1$. In particular:

\nero  all (absolute) de Rham cohomology groups $H^q(\Omega,\reals)$ vanish for $q\geq p$;

\nero  if  $p\leq \frac n2$ then  $H^q(\Sigma,\reals)=0$ for all $p\leq q\leq n-p$. 

\medskip

\acapo (The second statement follows by looking at the long-exact cohomology sequence of the pair $(\Omega,\Sigma)$).

In section \ref{tf} we will observe that these vanishing results hold under the assumption that $W^{[p]}_{\Omega}\geq 0$, which is sometimes weaker than assuming the non-negativity of the sectional curvature.  
\end{rem}


\subsection{Rigidity results}


From the equality case of the main theorem we see that  any geometric situation in which the equality holds gives rise to a rigidity result, that is, implies that $\Omega$ is a Euclidean ball. We study the following two situations:

\parte a The boundary $\Sigma$ carries a non-trivial special Killing form (in particular, if $\Sigma$ is isometric to a round sphere).

\parte b  The ambient manifold $\Omega$ carries a non trivial parallel $p$-form for some $p=1,\dots, n$ and $\Sigma$ is totally umbilical.

\medskip

We prove that, under suitable curvature conditions, $\Omega$ must be isometric to a Euclidean ball: see for example Theorem
\ref{boundarysphere}, Corollary \ref{ballrigidity} and Theorem \ref{rigidity2}. In particular, we observe that:

\medskip

\it If $\Omega$ has non-negative sectional curvature and admits a non-trivial parallel $p$-form for some $p=1,\dots,n$,  and $\Sigma$ is totally umbilical and has positive constant mean curvature, then $\Omega$ is isometric to a Euclidean ball. 

\medskip\rm

We refer to Section \ref{rigidity} for the definition of special Killing form and the precise statements of our results.
 
 
\subsection{Upper bounds} 


As a by-product of our calculations, we also observe some upper bounds of the first Hodge eigenvalue of the boundary $\Sigma$ when the ambient manifold $\Omega$ carries a non-trivial parallel $p$-form. In Theorem \ref{boundpi} we prove that, if $2\leq p\leq n-1$ and $H^p(\Sigma,\reals)=H^{n-p+1}(\Sigma,\reals)=0$, then:
$$
\lambda_{1,p}'(\Sigma)\leq \max\{p,n-p+1\}\dfrac{\int_{\Sigma}\norm S^2}{{\rm Vol}(\Sigma)},
$$
where $\norm{S}^2$ is the squared norm of the second fundamental form (if $\Sigma$ is minimal then the constant can be improved). Note the topological assumption; however we make no assumption on the curvatures of $\Omega$ and $\Sigma$. We will actually prove a slightly stronger inequality, which is sharp when $p=\frac{n+1}2$.
 
We refer to Theorem \ref{boundone} for a similar estimate in degree $p=1,n$: namely, if $H^1(\Sigma,\reals)=0$ and $\Omega$ supports a non-trivial parallel $1$-form, then
$$
 \lambda_1(\Sigma)\leq n\dfrac{\int_{\Sigma}\norm S^2}{{\rm Vol}(\Sigma)},
$$
where $\lambda_1(\Sigma)$ is the first positive eigenvalue of the Laplacian on functions.


\section{The Reilly formula for $p$-forms}\label{ReillyFormula}


Let $\Omega$ be an $(n+1)$-dimensional Riemannian domain with smooth boundary $\bd\Omega=\Sigma$ and $\omega$ a smooth differential $p$-form on $\Omega$. The Hodge Laplacian acts on $\omega$ and has the well-known expression:
$$
\Delta\omega=d\delta\omega+\delta d\omega,
$$
where $\delta=d^{\star}$ is the adjoint of the exterior derivative $d$ with respect to the canonical inner product of forms. Recall the Bochner formula:
\begin{equation}\label{bochner}
\Delta\omega=\nabla^{\star}\nabla\omega + W^{[p]}_{\Omega}(\omega),
\end{equation}
where $W^{[p]}_{\Omega}$, the \it curvature term, \rm is a self-adjoint endomorphism acting on $p$-forms. From the work of Gallot-Meyer (see \cite{gallot-meyer}) we know that, if the eigenvalues of the curvature operator are bounded below by $\gamma\in\reals$, then  
$$
\scal{W^{[p]}_{\Omega}(\omega)}{\omega}\geq p(n+1-p)\gamma\norm{\omega}^2.
$$
In particular, we observe that:

\nero If the curvature operator of $\Omega$ is non-negative then $W^{[p]}_{\Omega}\geq 0$ for all $p=1,\dots,n$.

\medskip

Recall that $W^{[1]}_{\Omega}$ is just the Ricci tensor of $\Omega$. The Reilly formula is obtained by integrating \eqref{bochner} on $\Omega$: for the precise statement, we need to introduce some 
additional notations. Given $N$, the inner unit normal vector field on $\Sigma$, the shape operator $S$ is defined as $S(X)=-\nabla_XN$
for all tangent vectors $X\in T\Sigma$; it admits a canonical extension acting on $p$-forms on $\Sigma$ and denoted by $S^{[p]}$. 
Explicitly, if $\omega$ is a $p$-form on $\Sigma$ one has:
\begin{eqnarray*}
S^{[p]}\omega(X_1,\dots,X_p)=\sum_{j=1}^p \omega(X_1,\dots,S(X_j),\dots,X_p),
\end{eqnarray*}
for tangent vectors $X_1,\dots,X_p\in T\Sigma$. Observe that $S^{[n]}$ is just multiplication by $nH$ where $H=\frac 1n {\rm tr}S$ is the mean curvature. We set, by convention, $S^{[0]}=0$.

\smallskip\acapo

It is clear from the definition that the eigenvalues of $S^{[p]}$ are precisely the $p$-curvatures of $\Sigma$: therefore we have immediately
\begin{equation}\label{espi}
\scal{S^{[p]}\omega}{\omega}\geq \sigma_p(\Sigma)\norm{\omega}^2
\end{equation}
at all points of $\Sigma$ and for all $p$-forms $\omega$, where $\sigma_p(\Sigma)$ is the lower bound of the $p$-curvatures defined previously.

If $J:\Sigma\to\Omega$ is the canonical inclusion, denote by $J^{\star}$ the restriction of differential forms to the boundary and let $i_N$ be the interior multiplication by $N$. At any point of the boundary one has $\norm{\starred J\omega}^2+\norm{i_N\omega}^2=\norm{\omega}^2$.
 
Recall the Hodge $\star$-operator of $\Omega$, mapping $p$-forms to $(n+1-p)$-forms and let $\star_{\Sigma}$ be the corresponding operator on forms of $\Sigma$. Here is the main computational tool of this paper. In what follows, we let $d^{\Sigma},\delta^{\Sigma},\Delta^{\Sigma}$ the differential, codifferential and Laplacian acting on forms of $\Sigma$, respectively. 

\begin{thm}\label{Reilly}\rm (Reilly formula)
\it \, In the above notation, let $\omega$ be a $p$-form on $\Omega$ with $p\geq 1$. Then:
\begin{eqnarray}\label{r}
\int_{\Omega}\norm{d\omega}^2+\norm{\delta\omega}^2=
\int_{\Omega}\norm{\nabla\omega}^2+\scal{W^{[p]}_{\Omega}(\omega)}{\omega}
+2\int_{\Sigma}\scal{i_N\omega}{\delta^{\Sigma}(J^{\star}\omega)}
+\int_{\Sigma}{\Cal B}(\omega,\omega),
\end{eqnarray}
where the boundary term has the following expression: 
$$
{\Cal B}(\omega,\omega)=\scal{S^{[p]}(J^{\star}\omega)}{J^{\star}\omega}
+\scal{S^{[n+1-p]}(J^{\star}\star\omega)}{J^{\star}\star\omega}
$$
or, equivalently:
$$
{\Cal B}(\omega,\omega)=\scal{S^{[p]}(J^{\star}\omega)}{J^{\star}\omega}
+nH\norm{i_N\omega}^2-\scal{S^{[p-1]}(i_N\omega)}{i_N\omega}.
$$
\end{thm}
 
The equivalence between the two expressions of the boundary term follows because $\starred J(\star\omega)$ is equal (up to sign) to $ \star_{\Sigma}i_N\omega$, and because of the identity $\star_{\Sigma}S^{[p]}+S^{[n-p]}\star_{\Sigma}=nH\star_{\Sigma}$ as operators on $p$-forms of $\Sigma$.

\medskip

For convenience, the proof of Theorem \ref{Reilly} will be given in the last section. Observe that, if $f$ is a smooth function, then the classical Reilly formula  is obtained by applying Theorem \ref{Reilly} to $\omega=df$ and by recalling that $W^{[1]}_{\Omega}={\rm Ric}_{\Omega}$. We write it for completeness:
\begin{equation}
\begin{aligned}\label{ReillyFunction}
\int_{\Omega}(\Delta f)^2 &=
\int_{\Omega}[ \norm{\nabla^2 f}^2 + {\rm Ric}_{\Omega}(\nabla f,\nabla f)]\\
& + \int_{\Sigma}\left[2\dfrac{\bd f}{\bd N}\Delta^{\Sigma}f+
\scal{S(\nabla^{\Sigma} f)}{\nabla^{\Sigma} f}+nH(\dfrac{\bd f}{\bd N})^2\right].
\end{aligned}
\end{equation}


\section{Lower bounds of eigenvalues of the Hodge Laplacian}



\subsection{Topological facts}\label{tf}


We now observe some topological consequences of $p$-convexity. These fact are not really new: besides the work of Wu, already cited in the introduction, we mention for example \cite{guerini-savo}, where the $p$-convexity assumption has been used to bound from below the spectrum of the Laplacian acting on $p$-forms of a manifold with boundary, for the absolute and relative conditions. 

First, we recall some well-known facts. Let $\Omega$ be a compact $(n+1)$-dimensional manifold with smooth boundary. The Hodge-de Rham theorem for manifolds with boundary asserts that $H^k(\Omega,\reals)$, the absolute cohomology space of $\Omega$ in degree $k$ with real coefficients, is isomorphic to the space of harmonic $k$-forms satisfying the absolute boundary conditions; equivalently, an absolute cohomology class is uniquely represented by a $k$-form on $\Omega$ satisfying the boundary problem:
\begin{equation}\label{harmonicfield}
\twosystem
{d\omega=\delta\omega=0\quad\text{on}\quad \Omega,}
{i_N\omega=0\quad\text{on}\quad \Sigma.}
\end{equation}
The relative cohomology space, denoted by $H^k_D(\Omega,\reals)$, is isomorphic to $H^{n+1-k}(\Omega,\reals)$ by Poincar\'e duality, which is induced at the level of forms by the Hodge $\star$-operator. Any class in $H^k_D(\Omega,\reals)$ is uniquely represented by a $k$-form satisfying
$$
\twosystem
{d\omega=\delta\omega=0\quad\text{on}\quad \Omega,}
{\starred J\omega=0\quad\text{on}\quad \Sigma.}
$$

\begin{thm}\label{topology}
Let $\Omega$ be a compact manifold with boundary $\Sigma$ and assume that $W^{[p]}_{\Omega}\geq 0$. 

\parte a If $\sigma_p(\Sigma)>0$ then $H^p(\Omega,\reals)=0$.

\parte b If $\sigma_p(\Sigma)\geq 0$ and $H^p(\Omega,\reals)\ne 0$ then $\Omega$ admits a non-trivial parallel $p$-form and $\sigma_p(x)=0$ for all $x\in\Sigma$. 
\end{thm}

\begin{proof} 
Let $\omega$ be a $p$-form, solution of the problem \eqref{harmonicfield}. We apply the Reilly formula of Theorem \ref{Reilly} to $\omega$ and obtain:
$$
0=\int_{\Omega}\norm{\nabla\omega}^2+\scal{W^{[p]}_{\Omega}(\omega)}{\omega}+
\int_{\Sigma}{\Cal B}(\omega,\omega).
$$
Now:
$$
\begin{aligned}
\int_{\Sigma}{\Cal B}(\omega,\omega)&=
\int_{\Sigma}\scal{S^{[p]}(\starred J\omega)}{\starred J\omega}\\
&\geq\int_{\Sigma}\sigma_p\norm{J^{\star}\omega}^2\\
&=\int_{\Sigma}\sigma_p\norm{\omega}^2
\end{aligned}
$$ 
by \eqref{espi}. The inner integral is non-negative, and is zero only when $\omega$ is parallel (hence of constant norm). It is now clear that $\sigma_p(\Sigma)>0$ forces $\omega$ to be identically zero on $\Sigma$, hence on $\Omega$, which proves a). For part b), there exists a non-trivial solution $\omega$ by assumption; this form has to be parallel and as 
$\int_{\Sigma} {\Cal B}(\omega,\omega)=0$ we see that $\sigma_p(x)=0$ on $\Sigma$. 
\end{proof}


\subsection{The main lower bound}\label{main}


As $\Delta^\Sigma$ commutes with $d^\Sigma$ and $\delta^\Sigma$, it preserves the space of exact and co-exact forms. Let $\lambda_{1,p}'(\Sigma)$ (resp. $\lambda_{1,p}''(\Sigma)$) be the first eigenvalue of the Hodge Laplacian when restricted to exact (resp. co-exact) $p$-forms of $\Sigma$, and let $\lambda_{1,p}(\Sigma)$ be the first positive eigenvalue of $\Delta^\Sigma$. The Hodge decomposition theorem implies that:
$$
\twosystem
{\lambda_{1,p}(\Sigma)=\min\{\lambda_{1,p}'(\Sigma),\lambda_{1,p}''(\Sigma)\}}
{\lambda_{1,p}''(\Sigma)=\lambda_{1,p+1}'(\Sigma)}
$$
the second identity being obtained by differentiating eigenforms. The Hodge duality gives $\lambda_{1,p}''(\Sigma)=
\lambda_{1,n-p}'(\Sigma)$. Given that, it is enough to estimate the eigenvalue $\lambda_{1,p}'(\Sigma)$ for $1\leq p\leq\frac{n+1}2$. 

\medskip

We now state our main theorem, under slightly weaker assumptions than those given in the introduction, namely, we only assume that $W^{[p]}_{\Omega}\geq 0$. As we already remarked, if the curvature operator is non-negative then $W^{[p]}_{\Omega}\geq 0$ for all $p$ (and of course $\rm{Ric}_{\Omega}\geq 0$). 
\begin{thm}\label{Hodge2}
Let $\Omega$ be a compact Riemannian manifold of dimension $n+1$ with smooth boundary $\Sigma$, and let $1\leq p\leq \frac{n+1}2$. Assume that $W^{[p]}_{\Omega}\geq 0$ and the $p$-curvatures of $\Sigma$ are bounded below by $\sigma_p(\Sigma)>0$. Then:
$$
\lambda_{1,p}'(\Sigma)\geq\sigma_p(\Sigma)\sigma_{n-p+1}(\Sigma).
$$
The equality holds for a Euclidean ball. If in addition $\Omega$ has non-negative Ricci curvature, then the equality holds if and only if $\Omega$ is a Euclidean ball. 
\end{thm} 

\begin{proof}
As $p\leq n-p+1$ by assumption, we have $\sigma_{n-p+1}(\Sigma)\geq \sigma_p(\Sigma)>0$. Let $\phi$ be a co-exact eigenform associated to $\lambda=\lambda_{1,p-1}''(\Sigma)$ and consider the exact $p$-eigenform $\omega=d^\Sigma\phi$ also associated to $\lambda$. The proof depends on the existence of a suitable extension of $\omega$. Precisely, by the results of Duff and Spencer (Theorem 2 p. 148 of \cite{duff}), there exists a $(p-1)$-form $\hat\phi$ on $\Omega$ such that $\delta d\hat\phi=0$ and $\starred J\hat\phi=\phi$ on $\Sigma$. Consider the $p$-form $\hat\omega=d\hat\phi$. Then $\hat\omega$ satisfies:
$$
\twosystem
{d\hat\omega=\delta\hat\omega=0\quad\text{on $\Omega$};}
{J^{\star}\hat\omega=\omega\quad\text{on $\Sigma$}.}
$$
We apply the Reilly formula (\ref{r}) to $\hat\omega$ and obtain:
\begin{equation}\label{proof1}
0=\int_{\Omega}\norm{\nabla\hat\omega}^2+\scal{W^{[p]}_{\Omega}\hat\omega}{\hat\omega}+2\int_{\Sigma}\scal{i_N\hat\omega}{\delta^{\Sigma}\omega}+
\int_{\Sigma}{\Cal B}(\hat\omega,\hat\omega).
\end{equation}
One has $\delta^{\Sigma}\omega=\delta^{\Sigma}d^{\Sigma}\phi=\lambda\phi$ and, by \eqref{espi}: 
$$
{\Cal B}(\hat\omega,\hat\omega)\geq \sigma_p(\Sigma)\norm{J^{\star}\hat\omega}^2
+\sigma_{n-p+1}(\Sigma)\norm{J^{\star}\star\hat\omega}^2.
$$
Now $\norm{J^{\star}\star\hat\omega}^2=\norm{i_N\hat\omega}^2$; since $\omega=d^\Sigma\phi$ is an eigenform associated to $\lambda$, 
one has $\int_{\Sigma}\norm{J^{\star}\hat\omega}^2=\int_{\Sigma}\norm{\omega}^2=\lambda\int_{\Sigma}\norm{\phi}^2$. 
The inner integral in \eqref{proof1} is non-negative because $W^{[p]}_{\Omega}\geq 0$ and we end-up with the following inequality:
\begin{equation}\label{main1}
0\geq \int_{\Sigma}2\lambda\scal{i_N\hat\omega}{\phi}
+\lambda\sigma_p(\Sigma)\norm{\phi}^2+\sigma_{n-p+1}(\Sigma)\norm{i_N\hat\omega}^2.
\end{equation}

\begin{rem}\label{parallel}
Note that, if the equality holds in \eqref{main1}, then $\hat\omega$ is parallel. 
\end{rem}
On the other hand, one has:
$$
\begin{aligned}
0&\leq\norm{i_N\hat\omega+\dfrac{\lambda}{\sigma_{n-p+1}(\Sigma)}\phi}^2\\
&=\norm{i_N\hat\omega}^2+\dfrac{2\lambda}{\sigma_{n-p+1}(\Sigma)}
\scal{i_N\hat\omega}{\phi}+\dfrac{\lambda^2}{\sigma_{n-p+1}(\Sigma)^2}
\norm{\phi}^2
\end{aligned}
$$
from which 
$$
2\lambda\scal{i_N\hat\omega}{\phi}\geq
-\sigma_{n-p+1}(\Sigma)\norm{i_N\hat\omega}^2-\dfrac{\lambda^2}
{\sigma_{n-p+1}(\Sigma)}\norm{\phi}^2.
$$
Note that, if the equality holds, then $i_N\hat\omega=-\frac{\lambda}{\sigma_{n-p+1}(\Sigma)}\phi$. Substituting in \eqref{main1} we end-up with the inequality
$$
0\geq \left(\lambda\sigma_p(\Sigma)-\dfrac{\lambda^2}
{\sigma_{n-p+1}(\Sigma)}\right)\int_{\Sigma}\norm{\phi}^2,
$$
and then, as $\phi$ is non-vanishing, we get $\lambda\geq\sigma_p(\Sigma)\sigma_{n-p+1}(\Sigma)$ as asserted. 

\smallskip

\acapo \bf The equality case. \rm It is well known that:
$$
\lambda_{1,p}'(\sphere n)=p(n-p+1),
$$
with multiplicity $\binom{n+1}p$ (see for example \cite{gallot-meyer}). As $\sigma_p(\Sigma)=p$ for all $p$, we see that we have equality for the Euclidean unit ball. 

Conversely, assume that the equality holds and ${\rm Ric}_\Omega\geq0$. Then $\hat\omega$ is parallel, 
$\lambda=\sigma_p(\Sigma)\sigma_{n-p+1}(\Sigma)$ and so:
\begin{equation}\label{equalityone}
i_N\hat\omega=-\sigma_p(\Sigma)\phi.
\end{equation}
Being parallel, $\hat\omega$ has constant norm, and we can assume that $\norm{\hat\omega}=1$ on $\Omega$. We apply the Stokes formula and observe that:
$$
\int_{\Omega}\norm{d\hat\phi}^2=\int_{\Omega}\scal{\hat\phi}{\delta d\hat\phi}-\int_{\Sigma}\scal{\starred J\hat\phi}
{i_Nd\hat\phi}.
$$
As $d\hat\phi=\hat\omega$ and $\delta d\hat\phi=0$, we obtain from \eqref{equalityone}:
\begin{equation}\label{equalitytwo}
{\rm vol}(\Omega)=\sigma_p(\Sigma)\int_{\Sigma}\norm{\phi}^2.
\end{equation}
On the other hand, $\phi$ is a co-exact eigenform of the Laplacian on $\Sigma$ associated to $\lambda$, so that:
$$
\begin{aligned}
\int_{\Sigma}\norm{\phi}^2&=\dfrac1{\lambda}\int_{\Sigma}\norm{d^\Sigma\phi}^2\\
&=\dfrac1{\lambda}\int_{\Sigma}\norm{\omega}^2\\
&=\dfrac1{\lambda}\left(\int_{\Sigma}\norm{\hat\omega}^2-\int_{\Sigma}\norm{i_N\hat\omega}^2\right)\\
&=\dfrac1{\lambda}{\rm vol}(\Sigma)-\dfrac{\sigma_p(\Sigma)^2}{\lambda}\int_{\Sigma}\norm{\phi}^2
\end{aligned}
$$
that is:
\begin{equation}\label{equalitythree}
\int_{\Sigma}\norm{\phi}^2=\dfrac{1}{\lambda+\sigma_p(\Sigma)^2}{\rm vol}(\Sigma).
\end{equation}
From \eqref{equalitytwo} and \eqref{equalitythree} we conclude:
$$
\dfrac{{\rm vol}(\Sigma)}{{\rm vol}(\Omega)}=\sigma_p(\Sigma)+\sigma_{n-p+1}(\Sigma).
$$
Now recall that, at any point $x\in\Sigma$ one has $\frac{\sigma_p(x)}{p}\leq \frac{\sigma_q(x)}{q}$ whenever $p\leq q$. As 
$\frac{\sigma_n(x)}{n}=H(x)$, the mean curvature, one has $\sigma_p(x)\leq pH(x)$ and then:
$$
\begin{aligned}
\sigma_p(\Sigma)+\sigma_{n-p+1}(\Sigma)&\leq \sigma_p(x)+\sigma_{n-p+1}(x)\\
&\leq pH(x)+(n-p+1)H(x)\\
&= (n+1)H(x).
\end{aligned}
$$
Therefore:
\begin{equation}\label{equalityfour}
H(x)\geq\dfrac{1}{n+1}\dfrac{\rm{vol}(\Sigma)}{\rm{vol}(\Omega)}
\end{equation}
for all $x\in\Sigma$. As the Ricci curvature of $\Omega$ is assumed non-negative, we can apply a result of Ros (Theorem 1 in \cite{ros})
which implies that 
$$
\inf_{x\in\Sigma}H(x)\leq \dfrac{1}{n+1}\dfrac{\rm{vol}(\Sigma)}{\rm{vol}(\Omega)}
$$
with equality if and only if $\Omega$ is a Euclidean ball. Hence \eqref{equalityfour} forces $\Omega$ to be isometric to a Euclidean ball, as asserted. 
\end{proof}


\section{Rigidity results}\label{rigidity}


We now state some rigidity results, which are consequences of the equality case in the main theorem.


\subsection{Manifolds with special Killing forms}\label{mwskf}


Here we study rigidity results for manifolds with boundary when the geometry of the boundary has special properties. 

In spin geometry, Hijazi and Montiel \cite{hijazi} (see also \cite{raulot}) assume the existence of Killing spinors on the boundary; they  prove that, under suitable curvature assumptions, a Killing spinor on the boundary extends to a parallel spinor on the domain. This imposes strong restrictions on the geometry of the domain: in particular, it has to be Ricci flat.

Here we address a similar problem for differential forms. As we will see, a natural assumption is the existence of a \it special Killing form \rm on $\Sigma$. This kind of differential forms were introduced by S. Tachibana and W. Yu \cite{tachibana} and compact, simply connected manifolds admitting such forms were classified by U. Semmelmann \cite{sem}. Noteworthy examples are given, besides the round spheres, by Sasakian manifolds.
\begin{defi}
Let $\Sigma$ be a $n$-dimensional Riemannian manifold. A special Killing $p$-form with number $c\in\mathbb{R}$ 
is a coclosed $p$-form $\omega$ on $\Sigma$ such that:
$$
\twosystem{\nabla^\Sigma_X\omega=\frac{1}{p+1}i_Xd^\Sigma\omega}{\nabla^\Sigma_X(d^\Sigma\omega)= -c(p+1)X^*\wedge\omega}
$$
for all vector fields $X$ with dual $1$-form $X^*$. 
\end{defi}
 
An easy computation shows that a non-trivial special Killing $p$-form is always a co-closed eigenform for the Hodge Laplacian acting on $p$-forms associated with the eigenvalue $c(p+1)(n-p)$ (in particular $c\geq0$). Thus, if $\Sigma$ carries a non-trivial special Killing $p$-form, then:
\begin{equation}\label{pf}
\lambda_{1,p}''(\Sigma)=\lambda_{1,p+1}'(\Sigma)\leq c(p+1)(n-p).
\end{equation}
As we already remarked, the standard round sphere $\sphere n$ carries special Killing forms: indeed
all co-closed eigenforms of the Hodge Laplacian associated with the eigenvalue $\lambda_{1,p}''(\sphere n)=(p+1)(n-p)$ are special Killing $p$-forms with number $c=1$.

From this discussion, it appears that Theorem \ref{Hodge2} applies when the boundary carries a special Killing form. Indeed, we have:
\begin{thm}\label{rig}
Let $\Omega$ be an $(n+1)$-dimensional compact manifold such that $W^{[p]}_{\Omega}\geq 0$ for some $1\leq p\leq \frac{n+1}2$. Assume that there exists a non-trivial special Killing $(p-1)$-form $\phi$ with number $c>0$ on the boundary $\Sigma$ and that 
$\sigma_p(\Sigma)\geq p\sqrt c$. Then:

\parte a $d^\Sigma\phi$ is the restriction of a parallel $p$-form on $\Omega$.

\parte b If in addition ${\rm Ric}_{\Omega}\geq 0$ then $\Omega$ is isometric to the  Euclidean ball of radius $\frac{1}{\sqrt c}$. 
\end{thm}

\begin{proof}
As $\sigma_p(\Sigma)\geq p\sqrt c$ and $1\leq p\leq \frac{n+1}2$ one has $\sigma_{n-p+1}(\Sigma)\geq (n-p+1)\sqrt c$. By the first part of the main theorem:
$$
\lambda_{1,p}'(\Sigma)\geq p(n-p+1)c.
$$
On the other hand  we have from (\ref{pf}) that 
$$\lambda_{1,p-1}''(\Sigma)=\lambda_{1,p}'(\Sigma)\leq p(n-p+1)c.$$ 
Then we have equality in Theorem \ref{Hodge}, and we see from the proof of this result that $d^{\Sigma}\phi$ must be the restriction to $\Sigma$ of a parallel $p$-form on $\Omega$. The statement \parte b is also immediate from our main theorem. 
\end{proof}

When the boundary is isometric to a round sphere we can remove the assumption on the Ricci curvature. Namely:
\begin{thm}\label{boundarysphere}  
Let $\Omega$ be an $(n+1)$-dimensional compact manifold such that $W^{[p]}_{\Omega}\geq 0$ for some $1\leq p\leq \frac{n+1}2$. Assume that $\Sigma$ is isometric to the unit round sphere $\sphere n$ and that $\sigma_p(\Sigma)\geq p$. Then $\Omega$ is isometric with the Euclidean unit ball.
\end{thm}

\begin{proof} 
As $\Sigma$ is isometric to the unit round sphere, we see that $\lambda_{1,p-1}''(\Sigma)=\lambda_{1,p}'(\Sigma)=p(n-p+1)$,
with multiplicity given by $\binom{n+1}{p}$. Then, from the proof of the main theorem, we see that any eigenform associated to 
$\lambda_{1,p}'(\Sigma)$ is the restriction to $\Sigma$ of a parallel $p$-form on $\Omega$. In particular ${\Cal P}^p(\Omega)$, 
the vector space of parallel $p$-forms on $\Omega$, has dimension at least $\binom{n+1}{p}$ because the 
restriction map $\starred J$ is injective on ${\Cal P}^p(\Omega)$. But, on a manifold of dimension $n+1$, $\binom{n+1}{p}$ is the maximal number of linearly independent parallel $p$-forms, and this maximum is achieved if and only if the manifold is flat. Then $\Omega$ is flat, and it is now straightforward to prove that $\Omega$ must be isometric to the Euclidean unit ball. 
\end{proof}

We observe the following corollary when the ambient manifold is locally conformally flat with non negative scalar curvature. This result could be compared to similar rigidity results obtained for spin manifolds in \cite{miao} or \cite{raulot}.
\begin{cor}\label{ballrigidity}
Let $\Omega$ be a $2m$-dimensional compact and connected Riemannian manifold with smooth boundary $\Sigma$. Assume that 
$\Omega$ is locally conformally flat and that its scalar curvature is non-negative. If the boundary is isometric to the round sphere $\sphere {2m-1}$ and satisfies $\sigma_m(\Sigma)\geq m$ then $\Omega$ is isometric to the Euclidean unit ball.
\end{cor}

\begin{proof} 
We take $p=m=\frac{n+1}2$ in the previous theorem, and so it is enough to show that $W^{[m]}\geq 0$.
From a result of Bourguignon \cite{bourguignon} (Proposition $8.6$), the curvature term for $m$-forms on locally conformally
flat manifolds is given by:
\begin{eqnarray*}
\scal{W^{[m]}(\omega)}{\omega}&=& \frac{m}{2(2m-1)}R\norm{\omega}^2
\end{eqnarray*}

\noindent where $R$ is the scalar curvature of $\Omega$ and then we have $W^{[m]}\geq 0$ as asserted.
\end{proof}


\subsection{Manifolds with parallel forms}\label{mwpf}


In this section we consider Riemannian manifolds admitting \it parallel forms: \rm by that we will mean that the manifold supports a non-trivial parallel form for some degree $p\ne 0, {\rm dim}\, M$. Noteworthy examples of such manifolds are given by Riemannian products and Kaehler manifolds. Other examples are given by simply connected, irreducible manifolds which are not symmetric spaces of rank $\geq 2$ and whose holonomy group is a proper subgroup of $SO_{n+1}$.  For more details on this subject, we refer to Chapter $10$ of \cite{besse}.   

We focus on the notion of  \it extrinsic sphere, \rm whose definition goes back to Nomizu, Yano and B.-Y. Chen. Here we consider the orientable, codimension one case:
\begin{defi} 
Given a Riemannian manifold $M$, we will say that the compact hypersurface $\Sigma$ of $M$ is an \rm extrinsic hypersphere 
\it if it is orientable, totally umbilical and has constant, non-zero  mean curvature $H$. 

If in addition $\Sigma$ is the boundary of a compact domain $\Omega$ and the mean curvature is positive, we will then say that $\Omega$ is an \rm extrinsic ball. \rm 
\end{defi}
\smallskip

It is of interest to know when extrinsic spheres are actually isometric to round spheres. Extrinsic spheres in Kaehler manifolds have been studied extensively, and have been classified (in the simply connected case) by Kawabata, Nemoto and Yamaguchi \cite{yam}. Other works of Okomura \cite{okomura} and Nemoto \cite{nemoto} study the case of extrinsic spheres in locally product Riemannian manifolds. Note that both of these manifolds carry parallel forms.

\smallskip

We first observe that if the ambient manifold $M$ carries a parallel $p$-form, then any extrinsic hypersphere $\Sigma$ 
immersed in $M$ carries non-trivial special Killing forms. We will denote by ${\Cal P}^{[p]}(M)$ the vector space of parallel 
$p$-forms on $M$ and by $\mathcal{K}^{p-1}(\Sigma)$ the vector space of special Killing $(p-1)$-forms with number $H^2$ on $\Sigma$, where $H$ is the (constant) mean curvature of $\Sigma$.  

\begin{prop} 
\parte a Let $M^{n+1}$ be a Riemannian manifold and $\Sigma^n$ an extrinsic hypersphere in $M$. If $\xi$ is a parallel 
$p$-form on $M$ then the normal part $i_N\xi$  is a special Killing $(p-1)$-form with number $H^2$. More generally, the linear map
$$
i_N:\mathcal{P}^p(M)\hookrightarrow\mathcal{K}^{p-1}(\Sigma)
$$
is injective.

\parte b Assume that $\Omega$ is an extrinsic ball with boundary $\Sigma$ such that $W^{[p]}_{\Omega}\geq 0$. Then 
$i_N:\mathcal{P}^p(\Omega)\hookrightarrow\mathcal{K}^{p-1}(\Sigma)$ is actually an isomorphism. 

\end{prop}

\begin{proof} Let $\xi$ be any $p$-form on $M$, and $\Sigma$ be any hypersurface of $M$. Then, for all vector fields $X$ on $\Sigma$ one has the formulae (see the proof of Lemma \ref{prep}):
$$
\twosystem
{\nabla_X^{\Sigma}(\starred J\xi)=\starred J(\nabla_X\xi)+S(X)^{\star}\wedge i_N\xi}
{\nabla_X^{\Sigma}(i_N\xi)=i_N\nabla_X\xi-i_{S(X)}\starred J\xi}
$$
If $\xi$ is parallel and $\Sigma$ is totally umbilical with constant mean curvature $H$ we have $S=H\cdot{\rm Id}$ and then:
$$
\twosystem
{\nabla_X^{\Sigma}(\starred J\xi)=HX^{\star}\wedge i_N\xi}
{\nabla_X^{\Sigma}(i_N\xi)=-Hi_X(\starred J\xi)}
$$
Recall that $H\ne 0$. We obtain easily (see Lemma \ref{prep} (ii)):
\begin{equation}\label{parallel1}
\twosystem
{\delta^\Sigma(J^{\star}\xi)=-(n-p+1)Hi_N\xi}
{d^{\Sigma}i_N\xi=-pH\starred J\xi.}
\end{equation}
Observe that the $p$-form $\starred J\xi$ is non trivial, because otherwise the first equation of \eqref{parallel1} would give 
$i_N\xi=0$ hence $\xi=0$ on $\Sigma$. But this is impossible because then, being parallel, $\xi$ would be zero on $M$. In a same way, $i_N\xi$ is non-trivial. Note that the first equation shows that $i_N\xi$ is co-exact and the second one that $\starred J\xi$ is exact. 
 
Let $\phi=i_N\xi$. Then $\starred J\xi=-\frac{1}{pH}d^{\Sigma}\phi$ and so: 
$$
 \nabla^{\Sigma}_X\phi=\frac1p i_Xd^{\Sigma}\phi, 
 \quad \nabla^{\Sigma}_Xd^{\Sigma}\phi=-pH^2X^{\star}\wedge\phi.
$$
Since $\phi=i_N\xi$ is co-closed, the above says precisely that $\phi$ is a special Killing $(p-1)$-form with number $H^2$ on $\Sigma$. Clearly $i_N$ is injective. Part a) now follows. 

We prove part b). Since the boundary is totally umbilical we have $\sigma_k(\Sigma)=kH$ for all $k=1,\dots,n$. From our main estimate we have then:
$$
\lambda_{1,p}'(\Sigma)\geq p(n-p+1)H^2.
$$
On the other hand, given a special Killing $(p-1)$-form $\phi$ on $\Sigma$, we know that $\phi$ is a co-closed eigenform associated to the eigenvalue $p(n-p+1)H^2$. Then: 
$$
\lambda_{1,p}'(\Sigma)\leq p(n-p+1)H^2.
$$
We have equality in Theorem \ref{Hodge2} and from its proof we see that $\phi$ is the normal part of a parallel $p$-form on $\Omega$. So the map $i_N$ is also surjective.
\end{proof}

We now can state another rigidity theorem, which is the main result of this section.
\begin{thm}\label{rigidity2}
Let $\Omega$ be an $(n+1)$-dimensional extrinsic ball having non-negative Ricci curvature and admitting a parallel $p$-form for some 
$p=1,\dots,n$. Suppose that at least one of the following conditions holds:

\item 1) $\Omega$ has non-negative sectional curvature;

\item 2) $W^{[p]}_{\Omega}\geq 0$;

\item 3) $\Omega$ is contractible;

\item 4) $H^p_D(\Omega,\reals)=0$.

\medskip

\acapo Then  $\Omega$ is isometric to a Euclidean ball. 
\end{thm}

\begin{proof} 
We will first prove the assertion under the condition 4), which will be shown to be the weakest of all. Let $\hat\omega$ be a parallel $p$-form on $\Omega$. From (\ref{parallel1}) we see that $d^{\Sigma}i_N\hat\omega=-pH\starred J\hat\omega$, so that, if $\phi=-\frac1{pH}i_N\hat\omega$ one has:
$$
i_N\hat\omega=-pH\phi\quad\text{and}\quad d^{\Sigma}\phi=\starred J\hat\omega.
$$
We also have $\Delta^{\Sigma}\phi=\lambda\phi$ with $\lambda=p(n-p+1)H^2$. As in the proof of the main theorem, we can extend 
$\phi$ to a $p$-form $\hat\phi$ on $\Omega$ satisfying:
$$\twosystem
{\delta d\hat\phi=0}
{\starred J\hat\phi=\phi}
$$
We claim that $d\hat\phi=\hat\omega$. In fact, the form $\hat h=d\hat\phi-\hat\omega$ satisfies $d\hat h=\delta\hat h=0$ on $\Omega$ and $J^{\star}\hat h=0$ on $\Sigma$, so it is a cohomology class in $H^p_D(\Omega,\reals)$. By our assumption we have indeed $\hat h=0$. 

At this point we proceed as in the proof of the main theorem (Theorem \ref {Hodge2}). Assuming that $\hat\omega$ has (constant) unit norm, we obtain:
$$
{\rm vol}(\Omega)=pH\int_{\Sigma}\norm{\phi}^2\quad\text{and}\quad 
\int_{\Sigma}\norm{\phi}^2=\dfrac{1}{\lambda+p^2H^2}{\rm vol}(\Sigma),
$$
which combined give 
$$
H=\dfrac{1}{n+1}\dfrac{{\rm vol}(\Sigma)}{{\rm vol}(\Omega)}.
$$
This, in turn,  implies that $\Omega$ is a Euclidean ball because ${\rm Ric}_{\Omega}\geq 0$.

It remains to show that any of the conditions 1), 2), 3) implies 4). Now 1) implies 4) by the result of Wu mentioned in the 
Introduction: in fact, as $\sigma_p(\Sigma)=pH>0$ for all $p$ we see that $H^p(\Omega, \reals)=H^p_D(\Omega, \reals)=0$ for all $p\ne 0, n+1$. We now assume 2). Again we have $\sigma_p(\Sigma)> 0$ for all $p$ and, by assumption, 
$W^{[n-p+1]}_{\Omega}=W^{[p]}_{\Omega}\geq 0$. Then by Theorem \ref{topology} we see that $H^{n-p+1}(\Omega,\reals)=0$ and by 
duality $H^p_D(\Omega,\reals)=0$ as well. So 2) implies 4). Finally 3) trivially implies 4). The proof is complete. 
\end{proof}

We conclude by observing the following immediate consequences.
\begin{cor}\label{extrinsicball}
\parte a Let $\Omega$ be an extrinsic ball in a manifold with positive sectional curvature. Then $\Omega$ admits no parallel forms in degree $p\ne 0,{\rm dim}\,\Omega$. 

\parte b Let $M$ be a  manifold with positive sectional curvature supporting a parallel form. Then $M$ has no embedded extrinsic hypersphere.
\end{cor}

\begin{proof}
\parte a If there is a parallel form, we see from Theorem \ref{rigidity2} that $\Omega$ must be isometric to a 
Euclidean ball: but this contradicts the assumption on the positivity of the sectional curvature of $\Omega$.

\parte b $M$ is compact, and the positivity of the sectional curvature implies that any embedded extrinsic hypersphere is the common boundary of two domains in $M$. On one of these, say $\Omega$, the mean curvature has to be positive, which implies that $\Omega$ is an extrinsic ball with a parallel form: this is impossible by a). 
\end{proof}

We finally remark that there exist extrinsic hyperspheres of manifolds with parallel forms, for example Kaehler manifolds, which are not isometric to round spheres. 

In fact, given any Sasakian manifold $(\Sigma,g)$, not isometric to a round sphere, consider the metric cone over $\Sigma$, which is the manifold $\hat\Sigma=\Sigma\times\reals^+$ with the metric $\hat g=r^2g+dr^2$. Then $\Sigma$ embeds isometrically into $\hat\Sigma$ as the hypersurface $r=1$, and it is easy to check that $\Sigma$ is totally umbilical with mean curvature of constant absolute value $1$. On the other hand, it is well-known that then $\hat\Sigma$ is a Kaehler manifold: if $\xi$ is the Killing $1$-form on $\Sigma$ defining the Sasakian structure, then the $2$ form $\hat\xi=rdr\wedge\xi+\frac12 r^2d\xi$ is the Kaehler form of $\hat\Sigma$.

More generally (see \cite{sem}, Lemma 4.5) any manifold supporting a special Killing form can be isometrically embedded as an 
extrinsic hypersphere in a manifold (the metric cone over $\Sigma)$ with a parallel form.


\section{Upper bounds for the Hodge Laplacian}


We finish by giving an upper bound of the Hodge-Laplace eigenvalues of a hypersurface $\Sigma^n$ of a Riemannian manifold $M^{n+1}$ supporting parallel $p$-forms.

We start from the case $p=1$ and assume that the hypersurface bounds a compact domain.

\begin{thm}\label{boundone} 
Assume that $\Sigma$ is connected and bounds a compact domain $\Omega$ carrying a parallel $1$-form. 

\parte a If $\Sigma$ is minimal then $H^1(\Sigma,\reals)\ne 0$. More generally,
$$
\dim (H^1(\Sigma,\reals))\geq\dim ({\Cal P}^1(\Omega)),
$$
where ${\Cal P}^1(\Omega)$ is the vector space of parallel $1$-forms on $\Omega$.

\parte b If $H^1(\Sigma,\reals)=0$ then:
$$
\lambda_1(\Sigma)\leq n \dfrac{\int_{\Sigma}\norm S^2}{{\rm Vol}(\Sigma)}.
$$
where $\lambda_1(\Sigma)$ is the first positive eigenvalue of the Laplacian on functions. 
\end{thm}

In higher degrees, we have the following estimate. We let $\norm{S}_p^2$ be the sum of the largest $p$ squared  principal curvatures, that is, if $\eta_1,\dots,\eta_n$ are the principal curvatures we define
$$
\norm{S}_p^2=\max_{i_1<\dots<i_p}\{\eta_{i_1}^2+\dots+\eta_{i_p}^2\}.
$$
Clearly $\norm{S}_p^2\leq\norm{S}_n^2=\norm{S}^2$ for all $p$.

\begin{thm}\label{boundpi} 
Let $\Sigma^n$ be a compact, connected immersed (not necessarily embedded) hypersurface of $M^{n+1}$, and assume that $M^{n+1}$ carries a parallel $p$-form for some $p=2,\dots,n-1$. If $H^p(\Sigma,\reals)=H^{n-p+1}(\Sigma,\reals)=0$ then
$$
\lambda_{1,p}'(\Sigma)\leq \alpha(p)\cdot\dfrac{\int_{\Sigma}{\norm S}_{\alpha(p)}^2}{{\rm Vol}(\Sigma)}.
$$
where $\alpha(p)=\max\{p,n-p+1\}$. 
If in addition $\Sigma$ is minimal then 
$$
\lambda_{1,p}'(\Sigma)\leq c(n,p)\cdot\dfrac{\int_{\Sigma}\norm S^2}{{\rm Vol}(\Sigma)}.
$$
where $c(n,p)=\dfrac1n\max\{p(n-p),(p-1)(n-p+1)\}$.
\end{thm}

\begin{rem}\rm When $p=\dfrac{n+1}2$ the estimate is sharp. In that case the upper bound becomes
$$
\lambda_{1,p}'(\Sigma)\leq p\cdot\dfrac{\int_{\Sigma}{\norm S}_{p}^2}{{\rm Vol}(\Sigma)},
$$
which is an equality when $\Sigma=\sphere{2p-1}$, isometrically immersed in $\real{2p}$: in fact, one has  ${\norm S}_{p}^2=p$, and  it is known that  $\lambda_{1,p}'(\sphere{2p-1})=p^2$. 
\end{rem}

Before giving the proofs we make some observations. If $\xi$ is a parallel $p$-form on $M^{n+1}$ then we have from Lemma \ref{prep}(ii):
\begin{equation}\label{commutation2}
\twosystem
{\delta^{\Sigma}(\starred J\xi)=S^{[p-1]}(i_N\xi)-nHi_N\xi}
{d^{\Sigma}i_N\xi=-S^{[p]}(\starred J\xi)}
\end{equation}
It follows that:
\begin{equation}\label{commutation3}
\twosystem
{\norm{\delta^{\Sigma}(\starred J\xi)}^2=\norm{S^{[n-p+1]}(\starred J\star\xi)}^2}
{\norm{\delta^{\Sigma}(\starred J\star\xi)}^2=\norm{S^{[p]}(\starred J\xi)}^2}
\end{equation}
In fact, one has $\starred J\star\xi=\pm \star_{\Sigma}(i_N\xi)$ hence:
$$
\delta^{\Sigma}(\starred J\star\xi)=\pm \delta^{\Sigma}(\star_{\Sigma}(i_N\xi))=
\pm\star_{\Sigma}d^{\Sigma}i_N\xi.
$$
Since $\star_{\Sigma}$ is norm-preserving, we obtain the second identity in \eqref{commutation3} by using the second identity in \eqref{commutation2}. The first formula in \eqref{commutation3} is obtained by duality. 

\medskip

If $\phi$ is any $p$-form on $\Sigma$, we observe that:
\begin{equation}\label{first}
\norm{S^{[p]}\phi}^2\leq p\norm{S}_p^2\norm{\phi}^2.
\end{equation}
and, if in addition $\Sigma$ is minimal:
\begin{equation}\label{second}
\norm{S^{[p]}\phi}^2\leq \dfrac{p(n-p)}{n}\norm{S}^2\norm{\phi}^2.
\end{equation}
For the proof of \eqref{first} and \eqref{second}, observe that the eigenvalues of $S^{[p]}$ are the $p$-curvatures of $\Sigma$, that is, all sums $\eta_{i_1}+\dots+\eta_{i_p}$.  Number the principal curvatures so that $\lambda^2=(\eta_1+\dots+\eta_p)^2$ is the largest eigenvalue of the endomorphism 
$(S^{[p]})^2$. Then, by the Cauchy-Schwarz inequality:
$$
\lambda^2=(\eta_1+\dots+\eta_p)^2\leq p(\eta_1^2+\dots+\eta_p^2),
$$
so $ \lambda^2\leq p\norm{S}_p^2$ and \eqref{first} follows. If $\Sigma$ is minimal one has also
$$
\lambda^2=(\eta_{p+1}+\dots+\eta_n)^2\leq (n-p)(\eta_{p+1}^2+\dots+\eta_n^2).
$$
and then:
$$
\frac{\lambda^2}{p}+\frac{\lambda^2}{n-p}\leq \norm{S}^2,
$$
from which we derive \eqref{second}. 

\medskip

\acapo\it Proof of Theorem \ref{boundone}. \rm Let $\xi$ be a parallel $1$-form. Then $d^{\Sigma}i_N\xi=-S^{[1]}(\starred J\xi)$. As $\xi$ is co-closed on $\Omega$ we have $\int_{\Sigma}i_N\xi=0$. If  $\starred J\xi=0$ on $\Sigma$ then $i_N\xi$ is constant on $\Sigma$, hence it vanishes. But this is impossible, because then $\xi=0$ on $\Sigma$, hence everywhere. So $\starred J\xi$ is not trivial and the restriction map $\starred J: {\Cal P}^1(\Omega)\to \Lambda^1(\Sigma)$ is one to one.

We now prove a). If $\Sigma$ is minimal and $\xi$ is parallel then \eqref{commutation2} shows $\starred J\xi$ is harmonic, hence $\starred J$ maps one to one into the subspace of harmonic $1$-forms and the assertion follows.

We prove b). The $1$-form $\starred J\xi$ is closed, hence exact. Then:
$$
\lambda_{1,1}'(\Sigma)\int_{\Sigma}\norm{\starred J\xi}^2\leq\int_{\Sigma}\norm{\delta^{\Sigma}(\starred J\xi)}^2\leq n\int_{\Sigma}\norm{S}^2\norm{i_N\xi}^2\\
$$
Recall that $\int_{\Sigma}i_N\xi=0$. Then:
$$
\begin{aligned}
\lambda_1(\Sigma)\int_{\Sigma}\norm{i_N\xi}^2&\leq\int_{\Sigma}\norm{d^{\Sigma}i_N\xi}^2\\
&\leq \int_{\Sigma}\norm{S^{[1]}\starred J\xi}^2\\
&\leq \int_{\Sigma}\norm{S}^2\norm{\starred J\xi}^2
\end{aligned}
$$
Summing the two inequalities and taking into account that $\lambda_{1,1}'(\Sigma)=\lambda_1(\Sigma)$ we get the assertion. 
\hfill$\square$
\medskip

\acapo\it Proof of Theorem \ref{boundpi}. \rm 
Let $\xi$ be a parallel $p$-form on $M$. We can assume that $\xi$ has constant unit norm. Now $\starred J\xi$ is closed because $H^p(\Sigma,\reals)=0$. So, by \eqref{commutation3} and \eqref{first}:
$$
\begin{aligned}
\lambda_{1,p}'(\Sigma)\int_{\Sigma}\norm{\starred J\xi}^2&\leq \int_{\Sigma}\norm{\delta^{\Sigma}(\starred J\xi)}^2\\
&=\int_{\Sigma}\norm{S^{[n-p+1]}(\starred J\star\xi)}^2\\
&\leq (n-p+1)\int_{\Sigma}\norm{S}_{n-p+1}^2\norm{\starred J\star\xi}^2.
\end{aligned}
$$
Now consider the parallel form $\star\xi$. Then $\starred J\star\xi$ is closed, hence exact by our assumptions. So:
$$
\begin{aligned}
\lambda_{1,n-p+1}'(\Sigma)\int_{\Sigma}\norm{\starred J\star\xi}^2&\leq \int_{\Sigma}\norm{\delta^{\Sigma}(\starred J\star\xi)}^2\\
&=\int_{\Sigma}\norm{S^{[p]}(\starred J\xi)}^2\\
&\leq p\int_{\Sigma}\norm{S}_p^2\norm{\starred J\xi}^2.
\end{aligned}
$$
By Poincar\'e duality $\lambda_{1,n-p+1}'(\Sigma)=\lambda_{1,p}'(\Sigma)$. Summing the two inequalities
and taking into account that $\norm{\starred J\xi}^2+\norm{\starred J\star\xi}^2=\norm{\xi}^2=1$ we get
$$
\lambda_{1,p}'(\Sigma)\cdot\Vol (\Sigma)\leq \alpha(p)\int_{\Sigma}\norm{S}_{\alpha(p)}^2.
$$
If $\Sigma$ is minimal we proceed as before using  inequality \eqref{second} instead of \eqref{first}.
\hfill$\square$
\medskip


\section{Proof of Theorem \ref{Reilly}}


The proof of the Reilly formula (\ref{r}) depends on the Stokes formula and certain commutation relations, stated in the following Lemma.

\begin{lemme}\label{prep} 
Let $\Omega$ be a domain with smooth boundary $\Sigma$, and let $\omega$ be a $p$-form on $\Omega$. 

\item (i) One has:
$$
\int_{\Omega}\norm{d\omega}^2+\norm{\delta\omega}^2=
\int_{\Omega}\scal{\omega}{\Delta\omega}+
\int_{\Sigma}\scal{i_N\omega}{\starred J(\delta\omega)}-
\scal{\starred J\omega}{i_Nd\omega}.
$$
\item (ii) As forms on $\Sigma$:
$$
\left\{\begin{aligned}
&\delta^\Sigma(\starred J\omega)=\starred J(\delta\omega)
+i_N\nabla_N\omega+S^{[p-1]}(i_N\omega)-nHi_N\omega\\
&d^{\Sigma}i_N\omega= -i_Nd\omega+\starred J(\nabla_N\omega)-S^{[p]}(\starred J\omega).
\end{aligned}\right.
$$
\end{lemme}

We can now prove the theorem. From the Bochner formula one gets:
\begin{equation}\label{one}
\int_{\Omega}\scal{\omega}{\Delta\omega}=
\int_{\Omega}\norm{\nabla\omega}^2+\scal{W^{[p]}_{\Omega}\omega}{\omega}
+\dfrac12\Delta\norm{\omega}^2.
\end{equation}
By the Green formula:
$$
\begin{aligned}
\int_{\Omega}\dfrac12\Delta\norm{\omega}^2&=
\int_{\Sigma}\scal{\nabla_N\omega}{\omega}\\
&=\int_{\Sigma}\scal{\starred J(\nabla_N\omega)}{\starred J\omega}+
\scal{i_N\nabla_N\omega}{i_N\omega}.
\end{aligned}
$$
Substituting in \eqref{one} and then in the Stokes formula of Lemma \ref{prep}(i)
one obtains:
\begin{equation}\label{two}
\begin{aligned}
\int_{\Omega}\norm{d\omega}^2+&\norm{\delta\omega}^2=
\int_{\Omega}\norm{\nabla\omega}^2+\scal{W^{[p]}_{\Omega}\omega}{\omega}\\
&+\int_{\Sigma}
\scal{i_N\omega}{\starred J(\delta\omega)}
-\scal{\starred J\omega}{i_Nd\omega}
+\scal{\starred J(\nabla\omega)}{\starred J\omega}
+\scal{i_N\nabla_N\omega}{i_N\omega}.
\end{aligned}
\end{equation}
From the first formula in Lemma \ref{prep}(ii):
\begin{equation}\label{three}
\begin{aligned}
\int_{\Sigma}\scal{i_N\omega}{\starred J(\delta\omega)}&=
\int_{\Sigma}
\scal{i_N\omega}{\delta^\Sigma(\starred J\omega)}+nH\norm{i_N\omega}^2\\
&-\int_{\Sigma}\left(\scal{i_N\nabla_N\omega}{i_N\omega}+\scal{S^{[p-1]}(i_N\omega)}{i_N\omega}\right),
\end{aligned}
\end{equation}
and from the second:
\begin{equation}\label{four}
\begin{aligned}
-\int_{\Sigma}\scal{\starred J\omega}{i_Nd\omega}&=
\int_{\Sigma}\scal{\starred J\omega}{d^{\Sigma}i_N\omega}-\scal{\starred J(\nabla_N\omega)}{\starred J\omega}+\scal{S^{[p]}(\starred J\omega)}{\starred J\omega}\\
&=\int_{\Sigma}\scal{\delta^{\Sigma}(\starred J\omega)}{i_N\omega}-\scal{\starred J(\nabla_N\omega)}{\starred J\omega}+\scal{S^{[p]}(\starred J\omega)}{\starred J\omega}.
\end{aligned}
\end{equation}
Substituting (\ref{three}) and (\ref{four}) in (\ref{two}) we finally get the statement of the theorem. 

\medskip

\noindent {\it Proof of Lemma \ref{prep}}. \rm The formula in (i) is a direct consequence of Stokes formula: if $\omega$ is a $(p-1)$-form and $\phi$ is a $p$-form, then:
$$
\int_{\Omega}\scal{d\omega}{\phi}=
\int_{\Omega}\scal{\omega}{\delta\phi}-\int_{\Sigma}
\scal{\starred J\omega}{i_N\phi}.
$$
We prove the second  formula in (ii). For all vector fields $X$ on $\Sigma$ one has the formulae:
\begin{equation}\label{derivative}
\twosystem
{\nabla_X^{\Sigma}(\starred J\omega)=\starred J(\nabla_X\omega)+S(X)^{\star}\wedge i_N\omega}
{\nabla_X^{\Sigma}(i_N\omega)=i_N\nabla_X\omega-i_{S(X)}\starred J\omega}
\end{equation}
which can be verified by a direct application of the Gauss formula $\nabla_XY=\nabla^{\Sigma}_XY+\scal{S(X)}{Y}N$. Fix $x\in\Sigma$ and let $(e_1,\dots,e_n)$ be an orthonormal basis of $T_x\Sigma$, so that $(e_1,\dots,e_n,N)$ is an orthonormal basis of $T_x\Omega$. Then, by definition:
$$
\twosystem
{\delta^{\Sigma}\starred J\omega=-\sum_{i=1}^ni_{e_i}
\nabla_{e_i}^{\Sigma}(\starred J\omega)}
{d^{\Sigma}i_N\omega=\sum_{i=1}^ne_i^{\star}\wedge\nabla_{e_i}^{\Sigma}(i_N\omega)}
$$
and using \eqref{derivative} one verifies the formulae in (ii). The details are straightforward, and we omit them.
\hfill$\square$



\vspace{0.8cm}     
Authors addresses:     
\nopagebreak     
\vspace{5mm}\\     
\parskip0ex     
\vtop{\hsize=6cm\noindent\obeylines}     
\vtop{     
\hsize=8cm\noindent     
\obeylines     
Simon Raulot,
Laboratoire de Math\'ematiques R. Salem
UMR $6085$ CNRS-Universit\'e de Rouen
Avenue de l'Universit\'e, BP.$12$
Technop\^ole du Madrillet
$76801$ Saint-\'Etienne-du-Rouvray, France}     
     
\vspace{0.5cm}     
     
E-Mail:     
{\tt simon.raulot@univ-rouen.fr }  

\vtop{\hsize=6cm\noindent\obeylines}     
\vtop{     
\hsize=9cm\noindent     
\obeylines     
Alessandro Savo,
Dipartimento di Metodi e Modelli Matematici 
Sapienza Universita' di Roma
Via Antonio Scarpa 16, 00161 Roma, Italy         
}     
     
\vspace{0.5cm}     
     
E-Mail:     
{\tt savo@dmmm.uniroma1.it  } 



\begin{thebibliography}{1}

\bibitem{besse} A.L. Besse, \emph{Einstein manifolds},  {Springer-Verlag, New-York} (1987).

\bibitem{bourguignon}    
J.-P. Bourguignon, \emph{Les vari\'et\'es de dimension $4$ \`a signature non nulle dont la courbure est harmonique sont d'Einstein},    
{Inv. Math.} \textbf{63} (1981), 263--286.

\bibitem{choi-wang}    
H.I. Choi and A.N. Wang, \emph{A first eigenvalue estimate for minimal hypersurfaces},    
{J. Diff. Geo.} \textbf{18} (1983), 559--562. 

\bibitem{duff}
G.F.D. Duff and D.C. Spencer, \emph{Harmonic tensors on Riemannian manifolds with boundary},
{Ann. Math.} \textbf{57}, 127-156

\bibitem{gallot-meyer}    
S. Gallot and D. Meyer, \emph{Op\'erateur de courbure et laplacien des formes diff\'erentielles d'une vari\'et\'e riemannienne},
 {J. Math. Pures. Appl.} \textbf{54} (1975), 259--284.

\bibitem{guerini-savo}    
P. Guerini and A. Savo, \emph{Eigenvalue and gap estimates for the Laplacian acting on $p$-forms},    
{Trans. Amer. Math. Soc.}\textbf{ 356}  (2004), 319-344.

\bibitem{hijazi}    
O. Hijazi and S. Montiel, \emph{Extrinsic Killing Spinors},    
{Math. Z.} \textbf{243} (2003), 337--347.

\bibitem{hijazi-montiel-zhang}    
O. Hijazi, S. Montiel and X. Zhang, \emph{Dirac operator on embedded hypersurfaces},    
{Math. Res. Lett.} \textbf{8} (2001), 195--208.

\bibitem{yam} 
N. Kawabata, H. Nemoto and S. Yamaguchi, \emph{Extrinsic spheres in Kaehler manifolds}.
{Michigan Math. J.} \textbf{31} (1984), 15-19.

\bibitem{miao}    
P. Miao, \emph{Positive mass theorem on manifolds admitting corners along a hypersurface},    
{Adv. Theor. Math. Phys.} \textbf{6} (2003), 1163--1182.

\bibitem{nemoto}
H. Nemoto, \emph{Extrinsic spheres in a locally product Riemannian manifold},    
{Tensor N.S.} \textbf{40} (1983), 159--162. 

\bibitem{okomura}    
M. Okumura, \emph{Totally umbilical hypersurfaces of a locally product Riemannian manifold},    
{Kodai. Math. Sem. Rep.} \textbf{19} (1967), 35--42. 

\bibitem{petersen}
P. Petersen, \emph{Riemannian Geometry}, {Graduate Texts in Mathematics, Springer-Verlag, New-York} (1998).

\bibitem{raulot}    
S. Raulot, \emph{Rigidity of compact Riemannian spin manifolds with boundary},    
{Lett. Math. Phys.} \textbf{86} (2008), 177--192.

\bibitem{reilly}    
R.C. Reilly, \emph{Application of the {H}essian operator in a {R}iemannian Manifold},    
{Indiana Univ. Math. J.} \textbf{26} (1977), 459--472.

\bibitem{ros}    
A. Ros, \emph{Compact Hypersurfaces with constant higher order mean curvatures},    
{Revista Mathem\'atica Iberoamericana} \textbf{3} (1987), 447--453.

\bibitem{schwarz}    
G. Schwarz, \emph{Hodge {D}ecomposition-{A} method for solving boundary value problems},    
{Lecture Notes in Mathematics, Springer} (1995).

\bibitem{sem}    
U. Semmelmann, \emph{Conformal Killing forms on Riemannian manifolds}.
{Math. Z.} \textbf{245} (2003) no 3, 503--527.

\bibitem{tachibana}    
S. Tachibana and W.N. Yu, \emph{On a Riemannian space admitting more than one Sasakian structures}.
{Tohoku Math. J. (2)} \textbf{22} (1970), 536--540.

\bibitem{wu}    
H. Wu, \emph{Manifolds of partially positive curvature},    
{Indiana Univ. Math. J.} \textbf{36} (1987), 525--548.

\bibitem{xia}    
C. Xia, \emph{Rigidity of compact manifolds with boundary and nonnegative {R}icci curvature},    
{Proc. Amer. Math. Soc.} \textbf{125} (1997) no 6, 1801--1806.


\end{thebibliography}
\end{document}